# European and Asian Greeks for Exponential Lévy Processes

Anselm Hudde[1] · Ludger Rüschendorf[2]

Received: 11 February 2022 / Revised: 13 January 2023 / Accepted: 16 January 2023

**Abstract**
In this paper we give easy-to-implement closed-form expressions for European and Asian Greeks for general $L^2$-payoff functions and underlying assets in an exponential Lévy process model with nonvanishing Brownian motion part. The results are based on Hilbert space valued Malliavin Calculus and extend previous results from the literature. Numerical experiments suggest, that in the case of a continuous payoff function, a combination of Malliavin Monte Carlo Greeks and the finite difference method has a better convergence behavior, whereas in the case of discontinuous payoff functions, the Malliavin Monte Carlo method clearly is the superior method compared to the finite difference approach, for first- and second order Greeks. Reduction arguments from the literature based on measure change imply that the expressions for the Greeks in this paper also hold true for generalized Asian options in particular for fixed and floating strike Asian options.

## 1 Introduction

One of the classical applications of Malliavin Calculus is to find closed-form expressions of the so-called Greeks, which are partial derivatives of the estimated option prices with respect to certain parameters like the initial price, the volatility, and the risk-free interest rate. Let $X$ denote the underlying financial asset, let $\Phi$ be the payoff function of the underlying option, let $T$ be the maturity of the option and let $r$ denote the risk-free interest rate. Then the value $V$ of the option (or expected outcome) can be expressed as

$$V_E = E[e^{-rT}\Phi(X_T)] \tag{1.1}$$

for a European option and

✉ Anselm Hudde
anselm.hudde@fb3.fra-uas.de

Ludger Rüschendorf
ruschen@stochastik.uni-freiburg.de

1 Frankfurt University of Applied Sciences, Frankfurt, Germany

2 University of Freiburg, Freiburg im Breisgau, Germany

$$V_A = E\left[e^{-rT}\Phi\left(\frac{1}{T}\int_0^T X_t dt\right)\right] \tag{1.2}$$

for an Asian option. Common payoff functions for put or call options, respectively, are $\Phi_p(x) = \max\{0, (K-x)\} = (K-x)^+$ and $\Phi_c(x) = \max\{0, (x-K)\} = (x-K)^+$, where $K$ is the **exercise** or **strike price**. But it also includes exotic options like binary options, which provide a fixed payoff $p$ if the strike K is reached (see Hull (2022), p.239-243): $\Phi_b(x) = p\mathbb{1}_{\{x\geq K\}}$.

The problem of determining closed-form expressions or approximations for the pricing of European and Asian options resp. for generalized Asian options has found a lot of interest and various methods have been developed in the ample literature on this subject. Most of the results were given first in the context of the Black Scholes model, then further on extended to some classes of models driven by particular Lévy processes like the CGMY model or the general hyperbolic model see e.g. Dufresne (2005) and finally given in general form for some general classes of semimartingale models. For a detailed overview of this development, we refer to Vecer (2002) and Vecer and Xu (2004), see also Albrecher (2004) and Albeverio and Lütkebohmert (2005).

The main methods developed in this context are based on the inversion of the Laplace transform (extending Geman and Yor (1993)), its connection to the fast Fourier transform, on the analytical expansion method of Linetsky (2004) and the (integro-)differential equation method of Rogers and Shi (1995), Vecer and Xu (2004) and Vecer (2014). Based on a change of measure technique as in Shreve and Vecer (2000) it was shown in Vecer (2002) and Vecer and Xu (2004) that the path dependency in the formulation of the Asian option pricing problem can be simplified to the case of European options with modified payoff functions where the underlying asset is driven by a special semimartingale process leading to an integro-differential equation where the stock price is driven by a process with independent increments. This reduction is given for generalized Asian options of the form

$$\left(\alpha\int_0^T X_t.d\lambda(t) - K_1X_T - K_2\right)_+, \tag{1.3}$$

which includes for $K_1 = 0$ the fixed strike option and for $K_2 = 0$ the floating strike option. The averaging factor $\lambda$ of finite variation includes for $\lambda(t) = \frac{t}{T}$ the case of continuously sampled Asian options and for $\lambda(t) = \frac{1}{n}\lfloor\frac{nt}{T}\rfloor$ the case of discretely sampled Asian options as well as options of inception at $t_0 \geq 0$. These reduction results imply that the results on Greeks in our paper formulated for European and Asian options also hold true for generalized Asian options as in formula (1.3) including also Asian options at the inception $t_0 \geq 0$.

The main results in our paper are concerned with extending the Malliavin method to the calculation of the Greeks of European and Asian options for exponential Lévy models.

Let $\Delta = \frac{\partial}{\partial x}V$ be the Greek $\Delta$, i.e., the derivative of the option value with respect to the initial value $x$ of the underlying process $X$. Then, for an Asian option $\Delta$ can be expressed as

$$\Delta = \frac{\partial}{\partial x}E\left[e^{-rT}\Phi\left(\frac{1}{T}\int_0^T X_t dt\right)\right]. \tag{1.4}$$

A straightforward method for the numerical approximation of Greeks is the Monte Carlo finite difference method (see for example Hull (2022), p.472-473). An alternative method is the Monte Carlo Malliavin method, which has numerical advantages for discontinuous payoff functions. Malliavin Calculus is used to find a stochastic weight $\pi$ such that the derivative of the option value can be expressed as

$$\Delta = E\left[e^{-rT}\Phi\left(\frac{1}{T}\int_0^T X_t dt\right)\pi\right]. \tag{1.5}$$

The value of $\Delta$ is then computed by a Monte Carlo method. This method was introduced in Fournié et al. (1999), where the integration by parts formula of Malliavin Calculus is applied to derive closed formulas for Greeks in the Black Scholes model. In this model, the price process is given by

$$dX_t = rX_t dt + \sigma X_t dW_t, \quad X_0 = x \tag{1.6}$$

where $r \in \mathbb{R}$ is the risk-free interest rate, $\sigma \in (0, \infty)$ the volatility and $x \in (0, \infty)$ is the initial condition. An introduction to this topic is given in Montero and Kohatsu-Higa (2003). As demonstrated in Fournié et al. (1999), Bavouzet and Messaoud (2006), Davis and Johansson (2006) and Xu et al. (2014) this approach leads to better convergence behavior for discontinuous payoff functions (for example digital options) in comparison to the finite difference method.

There are various generalizations of the formulas of European Greeks to more general jump diffusions. Davis and Johansson (2006) generalize the Malliavin Calculus to jump diffusions of the form

$$dX_t = b(X_{t-})dt + \sigma(X_{t-})dW_t + \sum_{k=1}^{m} \alpha_k(X_{t-})(dN_t^{(k)} - \lambda_k dt), \quad X_0 = x \tag{1.7}$$

where $W$ is a Brownian motion and $N^{(1)}, \dots, N^{(m)}$ are Poisson processes.

Davis and Johansson (2006) assume a separability condition on the process, which states that there is a continuously differentiable function $f$ with bounded derivative in the first argument such that

$$X_t = f(X_t^c, X_t^d), \quad X_0^c = x \tag{1.8}$$

holds true. Davis and Johansson (2006) derive stochastic weights for European options.

On the other hand, Forster et al. (2009) do not make a separability assumption but rely instead on a Hörmander condition in order to ensure the existence of an ordinary (non-jump) diffusion between two jumps. Again, stochastic weights are given for European options.

Kawai and Takeuchi (2011) provide explicit formulas for European Greeks where the asset price dynamics is described by gamma processes and Brownian motions time-changed by a gamma process.

The most general results on stochastic weights for European options and non-pure jump processes in the literature so far is given in Petrou (2008). In this article, stochastic weights of the Greeks for European options based on market models given by solutions of the Lévy stochastic differential equation

$$dX_t = b(t, X_t)dt + \sigma(t, X_t)dW_t + \int_{\mathbb{R}_0} \gamma(t, z, X_{t-})\widetilde{\mu}(dz, dt), \quad X_0 = x \tag{1.9}$$

are given. Here $\widetilde{\mu}$ is a compensated Poisson random measure and the coefficients are assumed to be continuously differentiable with bounded derivatives. The Greeks are provided for the general case (1.9), and for a general stochastic volatility model with jumps both in the underlying and the volatility. In both cases, the formulas are not given as closed expressions, in that they contain stochastic integrals, or even Skorohod integrals.

Our paper is based on Petrou (2008). It is also based on the Hilbert space approach to Malliavin Calculus in particular on the chain rule, the Skorohod integral and the integration by parts formula. We give for the case of an underlying exponential Lévy model explicit formulas for the stochastic weights for European and for Asian options. The advantage of the formulas provided in this paper is that they are closed, i.e., without stochastic integrals. The implementation of a Monte Carlo Simulation of the formulas in this article is therefore straightforward. The main theoretical contribution of this paper, however, is that we also provide formulas for Asian options in Section 3.2.

More recent literature on Malliavin Monte Carlo weights include stochastic volatility models. Yilmaz (2018) derives Malliavin Monte Carlo Greeks for European options under the assumption that the underlying asset and interest rate both evolve from a stochastic volatility model and a stochastic interest rate model, respectively. Benth et al. (2021) consider an infinite dimensional Heston stochastic volatility model and derive stochastic weights for forward contracts.

There are also several papers which are concerned with discretely-averaged Asian options: Saporito (2020) derives a Monte Carlo approximation for the price of path-dependent derivatives of options under a multiscale stochastic volatility model.

In the case of continuously monitored Asian options, only few closed-form expressions for stochastic weights are given in the literature.

Montero and Kohatsu-Higa (2003) determine stochastic weights for Asian options in the Black Scholes model using classical Malliavin Calculus. El-Khatib and Privault (2004) and Huehne (2005) give closed formulas for the Greeks when the underlying asset is represented by a pure jump Lévy-process. Bavouzet and Messaoud (2006) derive Malliavin Monte Carlo weights in the Lévy case for jump diffusions of the form

$$S_t = \beta + \sum_{i=1}^{J_t} c(T_i, \Delta_i, S_{T_i^-}) + \int_0^t b(u, S_u)du + \int_0^t \sigma(u, S_u)dW_u \tag{1.10}$$

where $\sigma$ is assumed to be linearly bounded, has bounded first and second derivatives, and it is assumed that there is an $\epsilon > 0$ such that $|\sigma(u,x)| \geq \epsilon$ for all $(u,x) \in (0,\infty)^2$.

Pflug and Thoma (2016) introduce a measure valued differentiation approach to calculate Greeks of exotic options, which include Asian options, for discrete-time Lévy processes.

There are also approaches to compute Greeks which do not use Malliavin calculus: Kirkby (2017) uses Fourier based methods for numerical European Greeks under general exponential Lévy processes. Fusai et al. (2011) and Kirkby (2016) use Fourier-based methods for the pricing the discretely-monitored Asian options. Kirkby (2016) also describes how to generalize these results to continuously averaged options by the Richardson extrapolation. Some Monte Carlo Methods which combine pathwise derivatives and likelihood ratio method estimators are presented in Glasserman and Liu (2010). There are also analytical formulas for European Greeks in exponential Lévy models as in Aguilar et al. (2020).

The structure of this article is as follows: In Section 2, we recollect some basic results on the real-valued Malliavin calculus for Hilbert space valued Malliavin calculus according to Solé et al. (2007) which allows a self-contained reading of this paper. Based on this method we derive in Section 3 formulas for the Greeks of European and Asian options. In Section 4 we compare the convergence properties of the Malliavin Greeks with the finite difference method numerically for a class of examples given by jump diffusion models.

## 2 Malliavin Calculus for Hilbert Space Valued Processes

Let $(\widetilde{X}_t)_{t\in[0,T]}$, $\widetilde{X}_0 = 0$ be a square-integrable Lévy process with nonvanishing Brownian motion part on a probability space $(\Omega, \mathcal{F}, P)$. Let $x \in (0, \infty)$ and let $X_t = x\exp(\widetilde{X}_t)$. Then there is a geometric Brownian motion with drift $\gamma$

$$X_t^{(1)} = \gamma t + \sigma W_t,\ \sigma \neq 0$$

on a probability space $(\Omega_W, \mathcal{F}_W, P_W)$ and a pure jump process $X^J$ on a probability space $(\Omega_J, \mathcal{F}_J, P_J)$ such that $(\Omega, \mathcal{F}, P) \simeq (\Omega_W, \mathcal{F}_W, P_W) \otimes (\Omega_J, \mathcal{F}_J, P_J)$ and such that under this identification

$$\widetilde{X}_t = \gamma t + \sigma W_t + X_t^J$$

holds true. We assume that $W$ is a standard Brownian motion. The process $X^J$ can furthermore be decomposed into the sum of a compound Poisson process $X^{(2)} = \alpha \sum_{i=1}^{N_t} Y_i$ and a square-integrable pure jump martingale $X^{(3)}$ that almost surely has a countable number of jumps on finite intervals, i.e.,

$$X_t^J = \alpha \sum_{i=1}^{N_t} Y_i + X_t^{(3)}. \tag{2.1}$$

The intensity of the Poisson process $N$ is denoted by $\lambda$.

By means of the natural identification $L^2(\Omega_W \times \Omega_J) \simeq L^2(\Omega_W, L^2(\Omega_J))$ every random variable $Y \in L^2(\Omega)$ can be regarded as a random variable on $\Omega_W$ with values in the Hilbert space $L^2(\Omega_J)$.

For applications to the frame of exponential Lévy processes as described above we next remind some basic notion of Hilbert space valued Malliavin calculus as developed in Solé et al. (2007) extending the real-valued Malliavin Calculus in Nualart (2006).

Let $\mathbb{D}^{1,2} \subseteq L^2(\Omega_W)$ denote the domain of the Malliavin derivative $D$ and $\delta$ its adjoint operator. The following definition in Solé et al. (2007) extends the notion of Malliavin derivative to the Hilbert space case.

**Definition 2.1** (Malliavin derivative). For $X = \sum_{i=1}^{n} F_i v_i$ in the tensor product $\mathbb{D}^{1,2} \otimes L^2(\Omega_J)$, i.e., $F_i \in \mathbb{D}^{1,2}$ and $v_i \in L^2(\Omega_J)$ for all $i = 1, \dots, n$, the Malliavin derivative $D$ is defined by

$$DX := \sum_{i=1}^{n} DF_i \otimes v_i. \tag{2.2}$$

D is a closable operator with domain in $L^2(\Omega_W, L^2(\Omega_J))$ and values in $L^2([0,T] \times \Omega_W, L^2(\Omega_J))$. The closure of its domain is denoted by $\mathbb{D}^{1,2}(L^2(\Omega_J))$.

We use the same symbol $D$ for the derivative operator on $\mathbb{D}^{1,2}$ and on $\mathbb{D}^{1,2}(L^2(\Omega_J))$, but it becomes clear from the context which operator is meant.

The product of two random variables $X = \sum_{i=1}^{n} F_i v_i \in \mathbb{D}^{1,2} \otimes L^2(\Omega_J)$ and $Y = \sum_{j=1}^{m} G_j w_j \in \mathbb{D}^{1,2} \otimes L^2(\Omega_J)$ is defined pointwise, i.e.,

$$XY(\omega)(\omega') := \sum_{1\le i\le n, 1\le j\le m} F_i(\omega) G_j(\omega) v_i(\omega') w_j(\omega'). \tag{2.3}$$

This definition extends naturally to the whole domain of $\mathbb{D}^{1,2}(L^2(\Omega_J))$ and leads to the following product formula.

**Proposition 2.2** (Product formula). *The Malliavin derivative on* $\mathbb{D}^{1,2}(L^2(\Omega_J))$ *satisfies the product rule: For* $X, Y \in \mathbb{D}^{1,2}(L^2(\Omega_J))$ *it holds* $XY \in \mathbb{D}^{1,2}(L^2(\Omega_J))$ *and*

$$D(XY) = XDY + YDX. \tag{2.4}$$

The chain rule can be generalized to differentiable functionals on $\mathbb{D}^{1,2}(L^2(\Omega_J))$ as established in Solé et al. (2007) for simple Lévy processes and in the general case in Petrou (2008).

**Theorem 2.3** (Chain rule). *Let* $\Phi : \mathbb{R}^m \to \mathbb{R}$ *be a continuously differentiable function with bounded partial derivatives and let* $X = (X^{(1)}, \dots, X^{(m)})$ *be a vector of random variables such that* $X^{(j)} \in \mathbb{D}^{1,2}(L^2(\Omega_J))$ *for all* $j = 1, \dots, m$. *Then it holds that*

$$D\Phi(X) = \sum_{i=1}^{m} \frac{\partial}{\partial x_i} \Phi(X) DX^{(i)}. \tag{2.5}$$

A basic notion of Malliavin Calculus is the divergence operator or Skorohod integral. It is characterized by a partial integration formula.

**Definition 2.4** The **divergence operator** or **Skorohod integral**

$$\delta : L^2([0,T] \times \Omega_W, L^2(\Omega_J)) \to L^2(\Omega_W, L^2(\Omega_J))$$

is defined as the adjoint operator of $D : L^2(\Omega_W, L^2(\Omega_J)) \to L^2([0,T] \times \Omega_W, L^2(\Omega_J))$. Its domain $\mathrm{dom}(\delta)$ is the set of $u \in L^2(\Omega_W, L^2(\Omega_J))$ such that there is a $c \in \mathbb{R}$ such that

$$E_{\Omega_W \times \Omega_J}\left[\int_0^T D_t X u_t dt\right] \le c\|X\|_{L^2(\Omega_W \times \Omega_J))} \tag{2.6}$$

for all $X \in \mathbb{D}^{1,2}$. The divergence operator is characterized by the integration-by-parts formula

$$E_{\Omega_W \times \Omega_J}\left[\int_0^T D_t X u_t dt\right] = E[X\delta(u)] \tag{2.7}$$

for all $X \in \mathbb{D}^{1,2}$.

From the Riesz representation theorem, it follows that inequality (2.6) describes the greatest possible domain on which $\delta$ can be defined. The dual operator is also sometimes referred to as the **adjoint operator** of $D$. Since the domain of $D$ is dense in $L^2(\Omega_W \times \Omega_J)$, this operator is well-defined. It turns out that the operator $\delta$ is unbounded and closed, and $\mathrm{dom}(\delta)$ is dense in $L^2(\Omega_W, L^2(\Omega_J))$.

The following proposition due to Petrou (2008) shows how a real-valued random variable can be factored out of the divergence. This will be a useful tool in order to find an explicit representation of the divergence in the following section.

**Proposition 2.5** *Let* $X \in \mathbb{D}^{1,2}(L^2(\Omega_J))$ *and* $u \in \mathrm{dom}(\delta)$ *such that*

$$Xu \in L^2([0,T] \times \Omega_W, L^2(\Omega_J)) \text{ and } X\delta(u) - \int_0^T (D_tX)u_t dt \in L^2(\Omega_W, L^2(\Omega_J)).$$

*Then* $Xu \in \mathrm{dom}(\delta)$ *and*

$$\delta(Xu) = X\delta(u) - \int_0^T (D_tX)u_t dt. \tag{2.8}$$

## 3 The Calculation of Greeks of European Options and Asian Options in Exponential Lévy Models

In this section, we derive the formulas for the calculation of European and Asian Greeks in exponential Lévy models. Since the European case is already treated in Petrou (2008), the main theoretical contribution of this this paper is Section 3.2.

Let the underlying asset be described by a square integrable exponential Lévy process $X_t = x\exp(\widetilde{X}_t)$ with nonvanishing Brownian motion part such that $X_0 = x$. We stress that that there are no further conditions on $X$. With the Lévy-Itô decomposition (Cont and Tankov (2003), Proposition 3.7), $X_t$ can be represented as

$$X_t = \exp\left(\gamma t + \sigma W_t + \alpha \sum_{i=1}^{N_t} Y_i + X_t^{(3)}\right). \tag{3.1}$$

Here, $\sum_{i=1}^{N_t} Y_i$ is a compound Poisson process with intensity $\lambda$ of the Poisson process $N_t$, and $X_t^{(3)}$ is a square-integrable pure jump martingale that almost surely has a countable number of jumps on finite intervals.

Since $\exp(\gamma t + \sigma W_t) \in \mathbb{D}^{1,2}$, $X_t$ is Malliavin derivable, and it follows that

$$DX_t = D\exp(\gamma t + \sigma W_t) \otimes \exp(X_t^J) = \sigma\exp(\gamma t + \sigma W_t)\mathbb{1}_{[0,t]} \otimes \exp(X_t^J) = \sigma X_t \mathbb{1}_{[0,t]}. \tag{3.2}$$

### 3.1 European Greeks in Exponential Lévy Models

As in the Brownian motion case, an extension of the **integration by parts formula** from Nualart (2006) is the main tool for the calculation of European Greeks. This extended formula is given in Petrou (2008).

**Theorem 3.1** (Integration by parts formula). *Let $G$ be a real-valued random variable, $F \in \mathbb{D}^{1,2}(L^2(\Omega_J))$ and let $u \in L^2([0,T] \times \Omega)$ such that $\int_0^T D_tFu_t dt \neq 0$ a.s. and $Gu(\int_0^T D_tFu_t dt)^{-1} \in \mathrm{dom}(\delta)$. For $\Phi \in C_{\mathbf{b}}^1(\mathbb{R})$ continuously differentiable with bounded derivative we have*

$$E[\Phi'(F)G] = E\left[\Phi(F)\delta\left(\frac{u_tG}{\int_0^T D_tFu_t dt}\right)\right]. \tag{3.3}$$

As consequence, the integration by parts formula leads to closed formulas for European Greeks. Some related formulas are given in the case of Lévy diffusions in Petrou (2008), and in the case of asset price dynamics described by gamma processes and Brownian motions time-changed by a gamma process in Kawai and Takeuchi (2011). In comparison to the diffusion case, the following formulas are much simplified and allow for an easy implementation.

**Theorem 3.2** *Let $\Phi \in L^2((0,\infty))$ such that $E[\Phi(X_T)] < \infty$, let $V_0 = E\left[e^{-rT}\Phi(X_T)\right]$ be the option value at time $T = 0$. Then the Greeks for European options of exponential Lévy processes are given by*

$$\Delta = \frac{\partial V_0}{\partial x} = \frac{e^{-rT}}{x\sigma T} E\left[\Phi\left(X_T\right) W_T\right] \tag{3.4}$$

$$\mathcal{V} = \frac{\partial V_0}{\partial \sigma} = e^{-rT} E\left[\Phi\left(X_T\right)\left(\frac{W_T^2}{\sigma T} - W_T - \frac{1}{\sigma}\right)\right] \tag{3.5}$$

$$\rho = \frac{\partial V_0}{\partial r} = Te^{-rT} E\left[\Phi\left(X_T\right)\left(\frac{W_T^2}{\sigma T} - 1\right)\right] \tag{3.6}$$

$$\Theta = -\frac{\partial V_0}{T} = -e^{-rT} E\left[\Phi\left(X_T\right)\left(\frac{W_T^2}{2T^2} + \mu \frac{W_T}{\sigma T} - \left(\frac{1}{2T} + r\right)\right)\right] \tag{3.7}$$

$$\Gamma = -\frac{\partial V_0}{\partial x^2} = -\frac{e^{-rT}}{x^2 \sigma T} E\left[\Phi\left(X_T\right)\left(\frac{W_T^2}{\sigma T} - W_T - \frac{1}{\sigma}\right)\right] \tag{3.8}$$

$$A = -\frac{\partial V_0}{\partial \alpha} = -\frac{e^{-rT}}{\sigma T} E\left[\Phi\left(X_T\right) W_T \sum_{i=1}^{N_T} Y_i\right]. \tag{3.9}$$

The proofs of equations (3.4)–(3.8) are similar to the proofs in the Black Scholes case, see Fournié et al. (1999). We will exemplarily give the proof for equation (3.4). The Greek $A$ (capital Alpha) is a Greek that does not exist in the Black Scholes model, its derivation is also given here.

***Proof of Equation (3.4)*** We consider $X_T$ as a function of the initial value $x$, $X_T : (0,\infty) \to (0,\infty)$, $x \mapsto X_T(x) = x \exp(\tilde{X}_T)$. In the first part of the proof, we pose the additional assumption that $\Phi \in C_b^1((0,\infty))$ is a continuously differentiable function with bounded derivative. Since $\frac{\partial}{\partial x} X_T = \frac{1}{x} X_T$, it holds that

$$\frac{\partial}{\partial x} E[e^{-rT}\Phi(X_T)] = \frac{e^{-rT}}{x} E[\Phi'(X_T)X_T].$$

An application of the integration by parts formula from Theorem 3.1 with $F = G = X_T$ and $u = DX_T$ yields

$$E[\Phi'(X_T)X_T] = E\left[\Phi(X_T)\delta\left(\frac{(DX_T)X_T}{\int_0^T (D_tX_T)^2dt}\right)\right] = \frac{1}{\sigma T}E[\Phi(X_T)\delta(\mathbb{1}_{[0,T]})] = \frac{1}{\sigma T}E[\Phi(X_T)W_T]. \tag{3.10}$$

For the second part of the proof, first assume that $\Phi \in L^2((0,\infty))$ is bounded, but with no further conditions. Let $(\varphi_n)_{n\in\mathbb{N}}$ be an approximation of unity, and for each $n \in \mathbb{N}$, let $\Phi^{(n)} = \Phi * \varphi_n$ be the convolution of $\Phi$ and $\varphi_n$. Then $\Phi^{(n)}$ converges to $\Phi$ in $L^2((0,\infty))$,

$$\Phi^{(n)} \xrightarrow{L^2((0,\infty))} \Phi. \tag{3.11}$$

Moreover, each $\Phi_n$ is continuously differentiable, i.e. $\Phi_n \in C^1_{\mathbf{b}}((0,\infty))$. From (3.11) it follows that there is a subsequence $(m) \subseteq \mathbb{N}$ such that $\Phi^{(m)}$ converges to $\Phi$ a.s. and in $L^2((0,\infty))$. From Young's inequality and the assumption that $\Phi$ is bounded it follows that

$$\|\Phi^{(m)}\|_{L^\infty} = \|\Phi * \varphi_m\|_{L^\infty} \leq \|\Phi\|_{L^\infty}\|\varphi_m\|_{L^1} = \|\Phi\|_{L^\infty} =: C < \infty,$$

i.e., $(\Phi^{(m)})_{m\in\mathbb{N}}$ is bounded uniformly by the constant $C$. From the bounded convergence theorem it follows that

$$\Phi^{(m)}(X_T) \to \Phi(X_T) \tag{3.12}$$

$P$-a.s. and in $L^2(P)$. Define $h_m(x) = E\left[\Phi^{(m)}\left(X_T(x)\right)\right]$, and define $h(x) = E\left[\Phi\left(X_T(x)\right)\right]$ for all $x \in (0,\infty)$. We conclude that for all $x \in (0,\infty)$, $h_m(x)$ converges to $h(x)$. From the first part of the proof it follows that

$$\frac{\partial}{\partial x}h_m(x) = \frac{1}{x\sigma T}E\left[\Phi^{(m)}(X_T)W_T\right]. \tag{3.13}$$

Now, let $g(x) = \frac{1}{x\sigma T}E\left[\Phi\left(X_T(x)\right)W_T\right]$. From equation (3.10) and from (3.9) it follows that

$$\left|\frac{\partial}{\partial x}h_m(x) - g(x)\right| \leq \frac{1}{x\sigma T}\Big(\underbrace{E\left[\Phi^{(m)}(X_T) - \Phi(X_T)\right]^2}_{:=\epsilon_m(x)}\Big)^{\frac{1}{2}}E\left[W_T^2\right]^{\frac{1}{2}} \xrightarrow{m\to\infty} 0 \tag{3.14}$$

for all $x \in (0,\infty)$. Since $\widetilde{X}$ is defined as a convolution with a Gaussian random variable, it has a density and therefore also $Y = \exp(\widetilde{X}_T)$ has a density $f_Y$ with respect to the Lebesgue measure $\lambda$. This implies

$$h(x) = E[\Phi(xY)] = \int_{-\infty}^{\infty}\Phi(xy)f_Y(y)dy = \frac{1}{x}\int_{-\infty}^{\infty}\Phi(z)f_Y\left(\frac{z}{x}\right)dz$$

for all $x \in (0,\infty)$. We conclude that $h$ is continuous. It follows that

$$E\left[\left(\Phi^{(m)}(xY) - \Phi(xY)\right)^2\right] = \epsilon_m(x)^2 \xrightarrow{m\to\infty} 0$$

uniformly on compacts $K \subseteq (0,\infty)$ and thus

$$\frac{\partial}{\partial x}h_m(x) \to g(x)$$

uniformly on compacts in $K \subseteq (0,\infty)$. A classical result from analysis finally implies, that $g$ is differentiable and that

$$h_n'(x) \xrightarrow{n\to\infty} h'(x) = g(x).$$

In the case $\Phi \in L^2((0,\infty))$ the sequence $\Phi_n = (-n \vee \Phi) \wedge n$ converges pointwise to $\Phi$ as $n \to \infty$. Furthermore, $\Phi_n(X_T) \to \Phi(X_T)$ in $L^2(\mu)$ where $\mu \sim X_T$, and also uniformly on compacts in $x$. The result follows by the previous approximation argument. This completes the proof of equation (3.4).

***Proof of Equation (3.9)*** As before in the first step we assume that $\Phi \in C^1_{\mathbf{b}}((0,\infty))$ is continuously differentiable with bounded derivative such that $E[\Phi(X_T)] < \infty$. The derivative of $X_T$ with respect to $\alpha$ is $\frac{\partial}{\partial\alpha}X_T = X_T\sum_{i=1}^{N_T} Y_i$. An application of the integration by parts formula from Theorem 3.1 with $u = \mathbb{1}_{[0,T]}$ yields

$$\begin{aligned}\frac{\partial}{\partial\alpha}E[\Phi(X_T)] &= E\left[\Phi'(X_T)X_T\sum_{i=1}^{N_T} Y_i\right] \\ &= E\left[\Phi(X_T)\delta\left(\frac{X_T\sum_{i=1}^{N_T} Y_i}{\sigma T X_T}\right)\right] \\ &= \frac{1}{\sigma T}E\left[\Phi(X_T)\delta\left(\sum_{i=1}^{N_T} Y_i\right)\right] \\ &= \frac{1}{\sigma T}E\left[\Phi(X_T)W_T\sum_{i=1}^{N_T} Y_i\right].\end{aligned}$$

The generalization to $\Phi \in L^2((0,\infty))$ such that $E[\Phi(X_T)] < \infty$ then follows analogously to the proof of equation (3.4).

### 3.2 Asian Greeks in Exponential Lévy Models

Asian options are options where the payoff is determined by the average price of the underlying asset. The option price for an Asian option with exercise time $T$ then is given by

$$V_0 = E\left[e^{-rT}\Phi\left(\frac{1}{T}\int_0^T X_t dt\right)\right],$$

where $\Phi$ is the payoff function. For all $n \in \mathbb{N}$ let

$$I_{(n)} = \int_0^T t^n X_t dt, \quad n \geq 0.$$

Before proceeding we present a useful lemma:

**Lemma 3.3** *The random variables $I_{(n)} = \int_0^T t^n X_t dt$ are in $\mathbb{D}^{1,2}(L^2(\Omega_J))$ and it holds:*

$$D_s I_{(n)} = D_s\int_0^T t^n X_t dt = \sigma\int_s^T t^n X_t dt \text{ for all } s \in [0,T], n \in \mathbb{N}, \tag{3.15}$$

$$\int_0^T D_s I_{(n)} ds = \sigma I_{(n+1)} \text{ for all } n \geq 0. \tag{3.16}$$

***Proof*** Equation (3.12) follows with an approximation of an increasing sequence of Riemann sums. Define the Riemann sums

$$S_{(n)}^k(X) = \frac{T}{2^k} \sum_{i=0}^{2^k-1} s_{i,k}^n X_{s_{i,k}}$$

where $s_{0,0} = 0$ and

$$s_{i,k+1} = \begin{cases} s_{\lfloor i/2 \rfloor, k} & \text{if } s_{\lfloor i/2 \rfloor, k} \leq \frac{Ti}{2^{k+1}} \\ \frac{Ti}{2^{k+1}} & \text{else} \end{cases}$$

for $i \leq 2^{k+1}$. The sequence $S_{(n)}^k(X)$ is clearly in $\mathbb{D}^{1,2}(L^2(\Omega_J))$ and converges to $I_{(n)}$ a.s. Since the sequence $\left(S_{(n)}^k(X)\right)_k$ is also monotonously increasing, it follows that $S_{(n)}^k(X)$ converges to $I_{(n)}$ in $L^2(\Omega)$. The derivative of $S_{(n)}^k(X)$ is

$$DS_{(n)}^k(X) = \frac{T}{2^k} \sum_{i=0}^{2^k-1} Ds_i^k X_{s_i} = \frac{\sigma T}{2^k} \sum_{i=0}^{2^k-1} s_i^k X_{s_i} \mathbb{1}_{[0,s_i]}.$$

$D_s S_{(n)}^k(X)$ is monotonously increasing and therefore also converges to $\int_s^T t^k X_t dt$ in $L^2(\Omega)$. It follows that $\int_0^T X_t dt \in \mathbb{D}^{1,2}(L^2(\Omega))$ and that

$$D_s \int_0^T t^n X_t dt = \sigma \int_s^T t^n X_t dt$$

holds true. (3.16) is a straightforward application of the integration by parts formula of the Riemann integral. □

As a consequence of Lemma 3.3, the following integration by parts formula is obtained, which will be the key for the calculation of Asian Greeks:

**Corollary 3.4** *Let Φ be a continuously differentiable payoff function with bounded derivatives and let F be a random variable such that $\frac{F}{I_{(1)}}$ is Skorohod integrable. Then it holds that*

$$E\left[\Phi'\left(\frac{I_{(0)}}{T}\right)F\right] = E\left[\Phi\left(\frac{I_{(0)}}{T}\right)\delta\left(\frac{TF}{\sigma I_{(1)}}\right)\right].$$

***Proof*** With Lemma 3.3 it follows that

$$\int_0^T D_s \Phi\left(\frac{I_{(0)}}{T}\right) ds = \int_0^T \Phi'\left(\frac{I_{(0)}}{T}\right) \frac{\sigma}{T}\left(\int_s^T X_t dt\right) ds = \frac{\sigma}{T}\Phi'\left(\frac{I_{(0)}}{T}\right) I_{(1)}.$$

The result then follows from the definition of the Skorohod integral in Definition 2.4.

**Theorem 3.5** (Asian Greeks) *Let* $\Phi \in L^2((0,\infty))$ *be such that* $E[\Phi(\frac{1}{T}\int_0^T X_t dt)] < \infty$. *Then the Greeks for Asian options are given by*

$$\mathcal{G} = e^{-rT} E\left[\Phi\left(\frac{I_{(0)}}{T}\right)\pi_{\mathcal{G}}\right], \quad \mathcal{G} \in \{\Delta, \mathcal{V}, \rho, \Theta, \Gamma, A\}$$

*where*

$$\pi_\Delta = \frac{1}{\sigma x}\left(-\sigma + W_T \frac{I_{(0)}}{I_{(1)}} + \sigma \frac{I_{(0)} I_{(2)}}{I_{(1)}^2}\right) \tag{3.17}$$

$$\pi_{\mathcal{V}} = \frac{1}{\sigma}\left(-(1+\sigma W_T) + \frac{W_T \int_0^T X_t W_t dt - \sigma \int_0^T t X_t W_t dt}{I_{(1)}} + \frac{\sigma(\int_0^T X_t W_t dt) I_{(2)}}{I_{(1)}^2}\right) \tag{3.18}$$

$$\pi_\rho = \left(\frac{W_T}{\sigma} - T\right) \tag{3.19}$$

$$\pi_\Theta = \left(r - \frac{1}{T} + \frac{\frac{1}{\sigma T} I_{(0)} W_T - \frac{1}{\sigma} W_T X_T + T X_T}{I_{(1)}} + \frac{\frac{1}{T} I_{(0)} I_{(2)} - I_{(2)} X_T}{I_{(1)}^2}\right) \tag{3.20}$$

$$\begin{aligned} \pi_\Gamma = \frac{1}{\sigma^2 x^2}\Bigg(& 2\sigma^2 - 4\sigma W_T \frac{I_{(0)}}{I_{(1)}} + ((W_T^2 - T) I_{(0)} - 4\sigma^2 I_{(2)}) \frac{I_{(0)}}{I_{(1)}^2} \\ & + \sigma(3 W_T I_{(2)} - \sigma I_{(3)}) \frac{I_{(0)}^2}{I_{(1)}^3} + 3\sigma^2 \frac{I_{(0)}^2 I_{(2)}^2}{I_{(1)}^4}\Bigg) \end{aligned} \tag{3.21}$$

$$\pi_A = \frac{1}{\alpha}\left(\frac{\frac{1}{\sigma} W_T \int_0^T X_t^{(2)} X_t dt - \int_0^T t X_t^{(2)} X_t dt}{I_{(1)}} + \frac{\int_0^T X_t^{(2)} X_t dt I_{(2)}}{I_{(1)}^2}\right). \tag{3.22}$$

We give the proofs of equations (3.17) and (3.22); the proofs of equations (3.18) - (3.21) can be found in the Appendix. Throughout the proofs, $\Phi \in C_{\mathbf{b}}^1((0,\infty))$ will be a continuously differentiable function with bounded derivative such that $E[\Phi(\frac{1}{T}\int_0^T X_t dt)] < \infty$ holds true. The generalization to general $\Phi \in L^2((0,\infty))$ then follows like in the proof of equation (3.4).

***Proof of Equation (3.17)*** The integration by parts formula in Corollary 3.4 yields

$$\Delta = \frac{\partial}{\partial x} E\left[e^{-rT}\Phi\left(\frac{I_{(0)}}{T}\right)\right] = e^{-rT} E\left[\Phi'\left(\frac{I_{(0)}}{T}\right)\frac{I_{(0)}}{xT}\right] = \frac{e^{-rT}}{\sigma x} E\left[\Phi\left(\frac{I_{(0)}}{T}\right)\delta\left(\frac{I_{(0)}}{I_{(1)}}\right)\right].$$

From Proposition 2.5 applied with $u = \mathbb{1}_{[0,T]}$, it follows that

$$\delta\left(\frac{I_{(0)}}{I_{(1)}}\right) = \frac{I_{(0)}}{I_{(1)}}\delta(\mathbb{1}_{[0,T]}) - \int_0^T D_s \frac{I_{(0)}}{I_{(1)}} ds. \tag{3.23}$$

We apply the chain rule in Theorem 2.3 to the second term of the right-hand side of equation (3.23), and obtain from Lemma 3.3:

$$\begin{aligned}\int_0^T D_s \frac{I_{(0)}}{I_{(1)}} ds &= \int_0^T \frac{I_{(1)}D_sI_{(0)} - I_{(0)}D_sI_{(1)}}{I_{(1)}^2} ds \\ &= \frac{\int_0^T D_sI_{(0)}ds}{I_{(1)}} - \frac{I_{(0)}\int_0^T D_sI_{(1)}ds}{I_{(1)}^2} = \sigma - \frac{\sigma I_{(0)}I_{(2)}}{I_{(1)}^2}.\end{aligned}$$

The divergence can therefore be expressed as

$$\delta\left(\frac{I_{(0)}}{I_{(1)}}\right) = W_T\frac{I_{(0)}}{I_{(1)}} - \sigma + \frac{\sigma I_{(0)}I_{(2)}}{I_{(1)}^2}.$$

This finally leads to a closed-form expression for $\Delta$:

$$\Delta = \frac{e^{-rT}}{\sigma x}E\left[\Phi\left(\frac{I_{(0)}}{T}\right)\left(-\sigma + W_T\frac{I_{(0)}}{I_{(1)}} + \sigma\frac{I_{(0)}I_{(2)}}{I_{(1)}^2}\right)\right].$$

***Proof of Equation (3.22)*** The derivative of $X_t$ with respect to $\alpha$ is

$$\frac{\partial}{\partial\alpha}X_t = X_t\frac{\partial}{\partial\alpha}\left(X_t^{(1)} + \alpha\sum_{i=1}^{N_t}Y_i + X_t^{(3)}\right) = X_t\sum_{i=1}^{N_t}Y_i = \frac{1}{\alpha}X_t^{(2)}X_t.$$

This gives us

$$\begin{aligned}\frac{\partial}{\partial\alpha}E\left[\Phi\left(\frac{I_{(0)}}{T}\right)\right] &= E\left[\Phi'\left(\frac{I_{(0)}}{T}\right)\frac{1}{\alpha T}\int_0^T X_t^{(2)}X_t dt\right] = \frac{1}{\alpha}E\left[\frac{\sigma}{T}\Phi'\left(\frac{I_{(0)}}{T}\right)I_{(1)}\left(\frac{\int_0^T X_t^{(2)}X_t dt}{\sigma I_{(1)}}\right)\right] \\ &= \frac{1}{\alpha}E\left[\int_0^T D_s\Phi\left(\frac{I_{(0)}}{T}\right)\left(\frac{\int_0^T X_t^{(2)}X_t dt}{\sigma I_{(1)}}\right)ds\right] \\ &= \frac{1}{\alpha}E\left[\Phi\left(\frac{I_{(0)}}{T}\right)\delta\left(\frac{\int_0^T X_t^{(2)}X_t dt}{\sigma I_{(1)}}\right)\right].\end{aligned}$$

The Skorohod integral can be calculated with Proposition 2.5, applied with $u = \mathbb{1}_{[0,T]}$:

$$\delta\left(\frac{\int_0^T X_t^{(2)}X_t dt}{\sigma I_{(1)}}\right) = \frac{\int_0^T X_t^{(2)}X_t dt}{\sigma I_{(1)}}W_T - \frac{1}{\sigma}\int_0^T D_s\frac{\int_0^T X_t^{(2)}X_t dt}{I_{(1)}}ds.$$

The second term can be rewritten as

$$\frac{1}{\sigma}\int_0^T D_s \frac{\int_0^T X_t^{(2)} X_t dt}{I_{(1)}} ds = \frac{1}{\sigma}\left( \frac{\int_0^T \int_s^T \sigma X_t^{(2)} X_t dt ds}{I_{(1)}} - \frac{\int_0^T X_t^{(2)} X_t dt \int_0^T D_s I_{(1)} ds}{I_{(1)}^2} \right)$$

$$= \frac{\int_0^T t X_t^{(2)} X_t dt}{I_{(1)}} - \frac{\int_0^T X_t^{(2)} X_t dt I_{(2)}}{I_{(1)}^2}.$$

As a consequence, it holds that

$$A = \frac{1}{\alpha} E\left[ \Phi\left(\frac{I_{(0)}}{T}\right) \left( \frac{\frac{1}{\sigma} W_T \int_0^T X_t^{(2)} X_t dt - \int_0^T t X_t^{(2)} X_t dt}{I_{(1)}} + \frac{\int_0^T X_t^{(2)} X_t dt I_{(2)}}{I_{(1)}^2} \right) \right].$$

***Remark 3.6*** a) The chain rule in Theorem 2.3 is the basis of the integration-by-parts formula leading to the formulas for the Greeks in Theorems 3.2 and 3.5. Since this rule also holds true for functions in the multidimensional case, the formulas for the Greeks can also be stated in a similar way for multidimensional Greeks.

b) As described in the introduction, the reduction results from the literature directly enable corresponding expressions for the Greeks also for generalized Asian options, as for example for discretely sampled Asian options at the inception time $t_0 \geq 0$. In particular it implies explicit formulas for fixed strike and for floating strike Asian options.

## 4 Numerical Example

In this section we use the R package Hudde (2021), which contains an implementation of the formulas (3.4)–(3.8) and (3.17)–(3.21) to compare the convergence properties of the Malliavin Greeks with the finite difference method in the Asian option case and for different payoff functions. In order to investigate the numerical properties of the Malliavin Monte Carlo Method in comparison with the finite difference method, we make use of a jump diffusion model allowing for more simple simulations. Consider the jump diffusion model

$$X_t = x \exp\left((r - \sigma^2/2)t + \sigma W_t\right) \exp\left( \alpha \sum_{i=1}^{N_t} Y_i \right), \tag{4.1}$$

where $N_t$ is a Poisson process, the $Y_i$ are such that $\sqrt{3} Y_i$ follow the Student t-distribution with 3 degrees of freedom (i.e., $sd(Y_i) = 1$), and $\alpha$ is a scale parameter for the size of the jumps. Note that posing of a no arbitrage condition is not done in model (4.1) as the formulas for the sensitivities (Greeks) also make sense for general drifts in the model. The parameters are $x = 100$, $r = -0.01$, $\sigma = 0.25$, $\alpha = 0.15$, and $\lambda = 1$. Since the distribution of $\exp(\int_0^T tX_t dt)$ is not known, we also need Monte Carlo simulation to calculate the option price, and hence for the finite difference method. In order to provide the best possible numeric differentiation for comparison, we use the R package Gilbert and Varadhan (2015) which uses Richardson's extrapolation, and which provides results that are more precise than the results obtained with the simple finite difference method. In all cases, we calculate the integrals with 252 discretization steps.

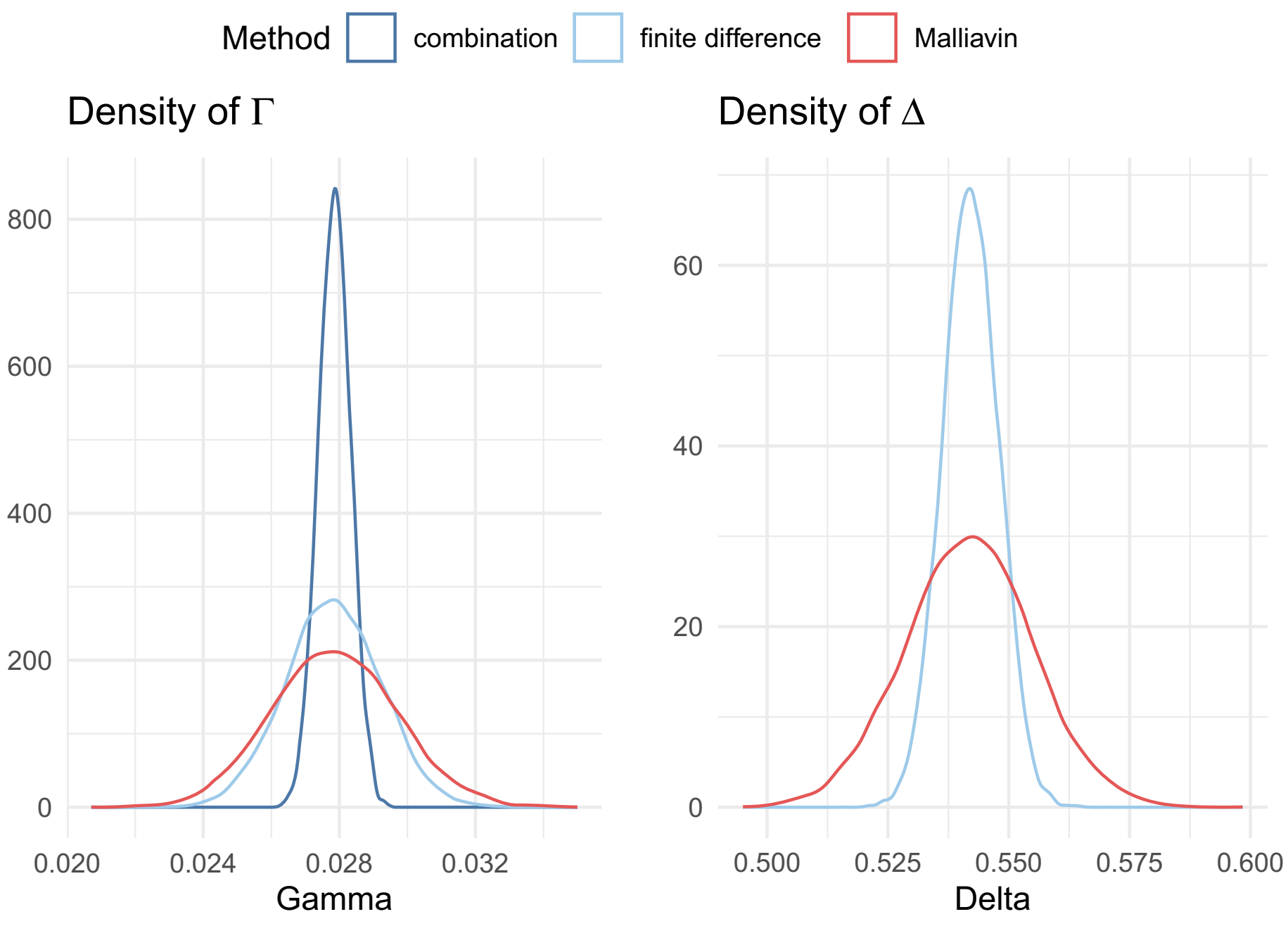

**Fig. 1** Densities of simulated $\Delta$ and $\Gamma$ for an Asian call option

## 4.1 Call Option in a Jump Diffusion Model

First, we consider the case of an Asian call option with continuous payoff function $\Phi_c(x) = \max\{0, (x - 100)\} = (x - 100)^+$. We run 10 000 simulations with 10 000 paths each. The results show that in this case, the finite difference method provides significantly better convergence properties for the first-order derivative Greeks, e.g., the Greek $\Delta$ (see Fig. 1), where the standard deviation of numeric $\Delta$ is 0.0058 compared to the standard deviation of the Malliavin $\Delta$ of 0.013 (see Table 1). In the case of second-order derivatives, the finite difference method still obtains better results, although the difference is much smaller (the standard deviation of the finite difference $\Gamma$ is 0.0014 in comparison to the standard deviation of the Malliavin $\Gamma$ of 0.0019). But the best result is obtained by a combination of both methods, i.e., by calculating $\Delta$ with the Malliavin Monte Carlo method, and then calculating its first-order derivative by the finite difference method. This results in a standard deviation of 0.00047 which decreases the standard deviation of the finite difference method by a factor of roughly 3.

**Table 1** Numeric results for an Asian call option

| | mean | sd | $q_{0.01}$ | $q_{0.99}$ |
|---|---|---|---|---|
| Malliavin $\Delta$ | 0.5418 | 0.01326 | 0.5108 | 0.5730 |
| finite difference $\Delta$ | 0.5421 | 0.005752 | 0.5288 | 0.5554 |
| Malliavin $\Gamma$ | 0.02786 | 0.001846 | 0.02368 | 0.03227 |
| finite difference $\Gamma$ | 0.02782 | 0.001405 | 0.02461 | 0.03110 |
| combi $\Gamma$ | 0.02786 | 0.0004719 | 0.02677 | 0.02895 |

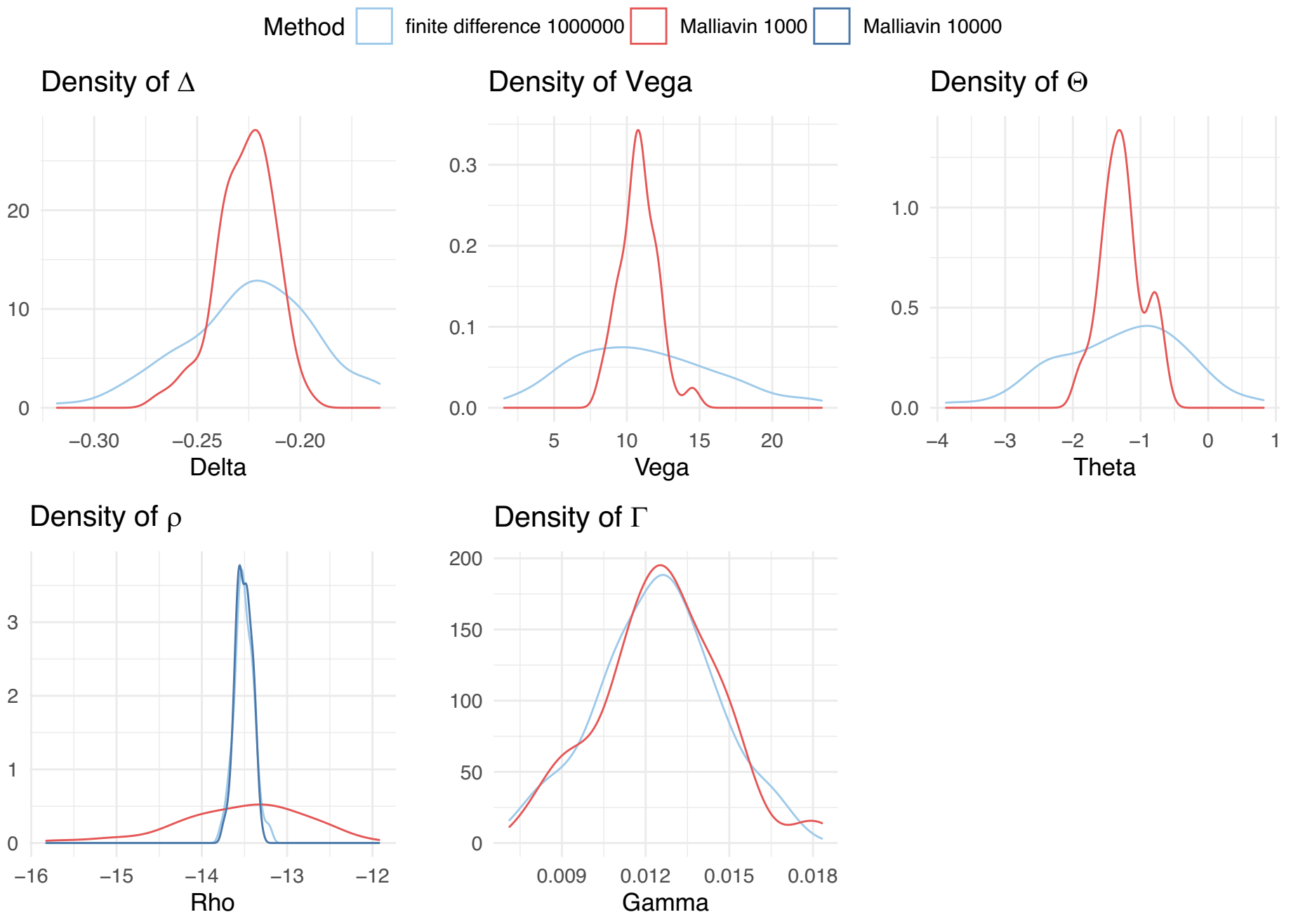

**Fig. 2** Densities of simulated Δ $\mathcal{V}$, Θ, $\rho$, and Γ for a digital Asian put option

## 4.2 Digital Option in a Jump Diffusion Model

Now we investigate an example of a discontinuous payoff function. Consider a digital put option with payoff function

$$\Phi(x) = \begin{cases} 10 & \text{if } x \leq 90 \\ 0 & \text{if } x > 90, \end{cases}$$

which results in a payoff of 10, if $\frac{1}{T}\int_0^T X_t dt \leq 90$. We run 100 simulations of the finite difference method with 1 000 000 paths, each, to compare the results with the Malliavin Monte Carlo method, where we run 100 simulations with 1 000 paths, each. For the Greek $\rho$, we also run 100 Simulations using the Malliavin method with 10 000 paths. The densities of the resulting distributions are plotted in Fig. 2, the aggregated results are also presented in Table 2. For all the plotted first-order Greeks except $\rho$, the variance of the results is still considerably smaller than in the case of the Malliavin Monte Carlo Method, although we used the squared number of simulated paths. For Γ, the results are nearly the same, and for $\rho$ the result with 10 000 comparable to the finite difference results with 1 000 000 paths. But still, the very large number of paths does not suffice to result in an acceptable accuracy.

We conclude that for the jump diffusion model introduced here, the Malliavin Monte Carlo provides a much better convergence behavior than the finite difference method.

**Table 2** Numeric results for a digital Asian put option

| | mean | sd | $q_{0.01}$ | $q_{0.99}$ |
|---|---|---|---|---|
| Malliavin $\Delta$ | -0.22663 | 0.013843 | -0.26396 | -0.20058 |
| finite difference $\Delta$ | -0.22444 | 0.031778 | -0.29782 | -0.16168 |
| Malliavin $\mathcal{V}$ | 10.834 | 1.2868 | 8.2257 | 14.427 |
| finite difference $\mathcal{V}$ | 11.217 | 4.8840 | 2.3248 | 22.612 |
| Malliavin $\Theta$ | -1.2714 | 0.30836 | -1.9196 | -0.67163 |
| finite difference $\Theta$ | -1.2876 | 0.94021 | -3.7941 | 0.63568 |
| Malliavin $\rho$ | -13.532 | 0.75671 | -15.660 | -12.182 |
| Malliavin $\rho_{10\,000}$ | -13.510 | 0.09239 | -13.734 | -13.340 |
| finite difference $\rho$ | -13.506 | 0.10995 | -13.762 | -13.220 |
| Malliavin $\Gamma$ | 0.012500 | 0.0021446 | 0.0078978 | 0.017840 |
| finite difference $\Gamma$ | 0.012387 | 0.0021373 | 0.0078785 | 0.016954 |

## 5 Summary

We have used a Hilbert-space valued Malliavin calculus to derive stochastic weights for sensitivities of options where the underlying asset is represented by a general exponential Lévy model with nonvanishing Brownian motion part. Results are given for European options and path-dependent Asian options, and for general $L^2$-payoff functions, in particular for common put, call and binary options. In the case of Asian options, this is a significant generalization of previous results in Montero and Kohatsu-Higa (2003).

These stochastic weights can be used for Monte Carlo simulations, which are easily implemented for arbitrary $L^2$-payoff functions and provide numerically good approximations of the Greeks.

We investigated numerically the quality of the Malliavin Monte Carlo method in comparison to the finite difference method for an easy-to-implement class of jump diffusion models. The simulations of the first-order Greeks with continuous payoff functions have a reasonable convergence behavior, and in this case the finite difference approach seems to be the better choice. For second-order Greeks a combination of the finite difference approach and the Malliavin Monte Carlo method for first-order Greeks is the superior method. Finally, for non-continuous payoff functions, like the case of binary options, the Malliavin Monte Carlo method provides a reasonable convergence behavior whereas the finite difference method is unstable even for extremely large numbers of paths. Thus, the stochastic weights introduced here provide a possibility to implement the computation of Greeks of digital Asian options for exponential Lévy processes, which seems to be very difficult with the finite difference method.

## Appendix

In the appendix we give the proofs of equations (3.18) to (3.21). Throughout the proofs, let $\Phi \in C^1_{\mathbf{b}}((0,\infty))$ be a continuously differentiable function with bounded derivative. The generalization to the payoff functions in $L^2((0,\infty))$ such that $E[\Phi(\frac{1}{T}\int_0^T X_t dt)] < \infty$ then follows analogously to the proof of equation (3.4).

***Proof of Equation (3.18)*** The derivative of $X_t$ with respect to $\sigma$ is

$$\frac{\partial}{\partial\sigma}X_t = X_t(W_t - \sigma t).$$

This gives us

$$\frac{\partial}{\partial\sigma}\frac{I_{(0)}}{T} = \frac{1}{T}\int_0^T X_t(W_t - \sigma t)dt = \frac{1}{T}\int_0^T X_t W_t dt - \frac{\sigma}{T}I_{(1)}.$$

With this representation, $\mathcal{V}$ can be written as

$$\begin{aligned}\mathcal{V} &= \frac{\partial}{\partial\sigma}E\left[e^{-rT}\Phi\left(\frac{I_{(0)}}{T}\right)\right]\\ &= e^{-rT}\underbrace{E\left[\frac{1}{T}\Phi'\left(\frac{I_{(0)}}{T}\right)\left(\int_0^T X_t W_t dt\right)\right]}_{(*)} - e^{-rT}\underbrace{E\left[\frac{\sigma}{T}\Phi'\left(\frac{I_{(0)}}{T}\right)I_{(1)}\right]}_{(**)}.\end{aligned}$$

The integration by parts formula in Corollary 3.4 then allows us to calculate

$$(*) = \frac{1}{\sigma}E\left[\Phi\left(\frac{I_{(0)}}{T}\right)\delta\left(\frac{\int_0^T X_t W_t dt}{I_{(1)}}\right)\right]$$

and

$$(**) = E\left[\Phi\left(\frac{I_{(0)}}{T}\right)\delta(\mathbb{1}_{(0,T]})\right] = E\left[\Phi\left(\frac{I_{(0)}}{T}\right)W_T\right].$$

We apply Proposition 2.5 with $u = \mathbb{1}_{[0,T]}$ and obtain

$$\begin{aligned}\delta\left(\mathbb{1}_{[0,T]}\frac{\int_0^T X_t W_t dt}{I_{(1)}}\right) &= \left(\frac{\int_0^T X_t W_t dt}{I_{(1)}}\right)\int_0^T \mathbb{1}_{[0,T]}(t)dW_t - \int_0^T D_s\left(\frac{\int_0^T X_t W_t dt}{I_{(1)}}\right)ds\\ &= \left(\frac{\int_0^T X_t W_t dt}{I_{(1)}}\right)W_T - \int_0^T D_s\left(\frac{\int_0^T X_t W_t dt}{I_{(1)}}\right)ds.\end{aligned}$$

With

$$D_s X_t W_t = X_t(\sigma W_t + 1)$$

for $s \le t$ it follows that

$$D_s\int_0^T X_t W_t dt = \sigma\int_s^T X_t W_t dt + \int_s^T X_t dt.$$

The chain rule in Theorem 2.3 yields

$$\int_0^T D_s\left(\frac{\int_0^T X_t W_t dt}{I_{(1)}}\right)ds$$
$$= \frac{I_{(1)}\Big(\sigma \int_0^T (\int_s^T X_t W_t dt)ds + \int_0^T (\int_s^T X_t dt)ds\Big) - (\int_0^T X_t W_t dt)\int_0^T D_s I_{(1)} ds}{I_{(1)}^2}$$
$$= \frac{\sigma \int_0^T tX_t W_t dt}{I_{(1)}} + 1 - \sigma\frac{(\int_0^T X_t W_t dt)I_{(2)}}{I_{(1)}^2}.$$

We can therefore write the divergence as

$$\delta\left(\mathbb{1}_{[0,T]}\frac{\int_0^T X_t W_t dt}{I_{(1)}}\right) = \left(\frac{\int_0^T X_t W_t dt}{I_{(1)}}\right)W_T - \frac{\sigma \int_0^T tX_t W_t dt}{I_{(1)}} - 1 + \sigma\frac{(\int_0^T X_t W_t dt)I_{(2)}}{I_{(1)}^2}.$$

This implies

$$(*) = \frac{1}{\sigma}E\left[\Phi\left(\frac{I_{(0)}}{T}\right)\left(-1 + \frac{W_T\int_0^T X_t W_t dt - \sigma\int_0^T tX_t W_t dt}{I_{(1)}} + \frac{\sigma(\int_0^T X_t W_t dt)I_{(2)}}{I_{(1)}^2}\right)\right].$$

It follows that

$$\mathcal{V} = \frac{e^{-rT}}{\sigma}E\Bigg[\Phi\left(\frac{I_{(0)}}{T}\right)\Bigg(-(1+\sigma W_T) + \frac{W_T\int_0^T X_t W_t dt - \sigma\int_0^T tX_t W_t dt}{I_{(1)}} + \frac{\sigma(\int_0^T X_t W_t dt)I_{(2)}}{I_{(1)}^2}\Bigg)\Bigg].$$

***Proof of Equation (3.19)*** We have $\frac{\partial}{\partial r}X_t = tX_t$ and therefore $\frac{\partial}{\partial r}I_{(0)} = I_{(1)}$. It follows that

$$\rho = \frac{\partial V_0}{\partial r} = \frac{\partial}{\partial r}E\left[e^{-rT}\Phi\left(\frac{I_{(0)}}{T}\right)\right] = -Te^{-rT}E\left[\Phi\left(\frac{I_{(0)}}{T}\right)\right] + e^{-rT}E\left[\Phi'\left(\frac{I_{(0)}}{T}\right)\frac{I_{(1)}}{T}\right].$$

The integration by parts formula in Corollary 3.4 yields

$$E\left[\Phi'\left(\frac{I_{(0)}}{T}\right)\frac{I_{(1)}}{T}\right] = \frac{1}{\sigma}E\left[\Phi\left(\frac{I_{(0)}}{T}\right)W_T\right],$$

and we conclude that

$$\rho = e^{-rT}E\left[\Phi\left(\frac{I_{(0)}}{T}\right)\left(\frac{W_T}{\sigma} - T\right)\right].$$

***Proof of Equation (3.20)*** A straightforward calculation gives us

$$\begin{aligned}\Theta = -\frac{\partial V_0}{\partial T} &= -\frac{\partial}{\partial T}E\left[e^{-rT}\Phi\left(\frac{I_{(0)}}{T}\right)\right] \\ &= re^{-rT}E\left[\Phi\left(\frac{I_{(0)}}{T}\right)\right] + e^{-rT}E\left[\Phi'\left(\frac{I_{(0)}}{T}\right)\left(\frac{I_{(0)}}{T^2} - \frac{X_T}{T}\right)\right] \\ &= rV_0 + \frac{x}{T}\Delta - e^{-rT}E\left[\Phi'\left(\frac{I_{(0)}}{T}\right)\frac{X_T}{T}\right].\end{aligned}$$

From the integration by parts formula in Corollary 3.4 it follows that

$$E\left[\Phi'\left(\frac{I_{(0)}}{T}\right)\frac{X_T}{T}\right] = E\left[\Phi\left(\frac{I_{(0)}}{T}\right)\delta\left(\frac{X_T}{\sigma I_{(1)}}\right)\right]. \tag{6.1}$$

To calculate the Skorohod integral in the expectation of the right-hand side of equation (6.1), we apply Proposition 2.5 with $u = \mathbb{1}_{[0,T]}$ and obtain

$$\delta\left(\frac{X_T}{\sigma I_{(1)}}\right) = \frac{X_T}{\sigma I_{(1)}}W_T - \int_0^T D_s\left(\frac{X_T}{\sigma I_{(1)}}\right)ds.$$

From the chain rule in Theorem 2.3 and from Lemma 3.3 it follows that

$$\int_0^T D_s\frac{X_T}{\sigma I_{(1)}}ds = \int_0^T \frac{\sigma^2 I_{(1)}X_T - \sigma\left(D_s I_{(1)}\right)X_T}{\sigma^2 I_{(1)}^2}ds = \frac{TX_T}{I_{(1)}} - \frac{I_{(2)}X_T}{I_{(1)}^2}.$$

Together, this leads to

$$\delta\left(\frac{X_T}{\sigma I_{(1)}}\right) = \frac{X_T}{\sigma I_{(1)}}W_T - \frac{TX_T}{I_{(1)}} + \frac{I_{(2)}X_T}{I_{(1)}^2} = \frac{\frac{1}{\sigma}W_T X_T - TX_T}{I_{(1)}} + \frac{I_{(2)}X_T}{I_{(1)}^2}.$$

We can finally write

$$\begin{aligned}\Theta &= -\frac{\partial}{\partial T}E\left[e^{-rT}\Phi\left(\frac{I_{(0)}}{T}\right)\right] \\ &= rV_0 + \frac{x}{T}\Delta - e^{-rT}E\left[\Phi'\left(\frac{I_{(0)}}{T}\right)\frac{X_T}{T}\right] \\ &= e^{-rT}E\left[\Phi\left(\frac{I_{(0)}}{T}\right)\left(r - \frac{1}{T} + \frac{\frac{1}{\sigma T}I_{(0)}W_T}{I_{(1)}} + \frac{\frac{1}{T}I_{(0)}I_{(2)}}{I_{(1)}^2} - \frac{(\frac{1}{\sigma}W_T - T)X_T}{I_{(1)}} - \frac{I_{(2)}X_T}{I_{(1)}^2}\right)\right] \\ &= e^{-rT}E\left[\Phi\left(\frac{I_{(0)}}{T}\right)\left(r - \frac{1}{T} + \frac{\frac{1}{\sigma T}I_{(0)}W_T - \frac{1}{\sigma}W_T X_T + TX_T}{I_{(1)}} + \frac{\frac{1}{T}I_{(0)}I_{(2)} - I_{(2)}X_T}{I_{(1)}^2}\right)\right]\end{aligned}$$

□

***Proof of Equation (3.21)*** To calculate $\Gamma$, we have to differentiate

$$\Delta = \frac{e^{-rT}}{\sigma x}E\left[\Phi\left(\frac{I_{(0)}}{T}\right)\left(-\sigma + W_T\frac{I_{(0)}}{I_{(1)}} + \sigma\frac{I_{(0)}I_{(2)}}{I_{(1)}^2}\right)\right]$$

with respect to the initial value $x$ again: This gives

$$\begin{aligned}\Gamma = \frac{\partial}{\partial x}\Delta &= -\frac{1}{x}\Delta + \frac{e^{-rT}}{\sigma x}E\left[\Phi'\left(\frac{I_{(0)}}{T}\right)\frac{I_{(0)}}{xT}\left(-\sigma + W_T\frac{I_{(0)}}{I_{(1)}} + \sigma\frac{I_{(0)}I_{(2)}}{I_{(1)}^2}\right)\right] \\ &+ \frac{e^{-rT}}{\sigma x}E\left[\Phi\left(\frac{I_{(0)}}{T}\right)\frac{\partial}{\partial x}\left(-\sigma + W_T\frac{I_{(0)}}{I_{(1)}} + \sigma\frac{I_{(0)}I_{(2)}}{I_{(1)}^2}\right)\right] \\ &= -\frac{1}{x}\Delta + \frac{e^{-rT}}{\sigma x}E\left[\Phi'\left(\frac{I_{(0)}}{T}\right)\frac{I_{(0)}}{xT}\left(-\sigma + W_T\frac{I_{(0)}}{I_{(1)}} + \sigma\frac{I_{(0)}I_{(2)}}{I_{(1)}^2}\right)\right]\end{aligned} \tag{6.2}$$

since $\frac{I_{(0)}}{I_{(1)}}$ and $\frac{I_{(0)}I_{(2)}}{I_{(1)}^2}$ are constant functions of $x$. An application of the integration by parts formula in Corollary 3.4 to the right-hand side of equation (6.2) yields

$$\Gamma = -\frac{1}{x}\Delta + \frac{e^{-rT}}{\sigma^2 x^2}E\left[\Phi\left(\frac{I_{(0)}}{T}\right)\delta\left(-\sigma\frac{I_{(0)}}{I_{(1)}} + W_T\frac{I_{(0)}^2}{I_{(1)}^2} + \sigma\frac{I_{(0)}^2 I_{(2)}}{I_{(1)}^3}\right)\right]. \tag{6.3}$$

An application of Proposition 2.5 with $u = \mathbb{1}_{[0,T]}$ to the Skorohod integral on the right-hand side of equation (6.3) with $u = \mathbb{1}_{[0,T]}$ yields

$$\begin{aligned}&\delta\left(-\sigma\frac{I_{(0)}}{I_{(1)}} + W_T\frac{I_{(0)}^2}{I_{(1)}^2} + \sigma\frac{I_{(0)}^2 I_{(2)}}{I_{(1)}^3}\right) \\ &= \left(-\sigma\frac{I_{(0)}}{I_{(1)}} + W_T\frac{I_{(0)}^2}{I_{(1)}^2} + \sigma\frac{I_{(0)}^2 I_{(2)}}{I_{(1)}^3}\right)\delta(\mathbb{1}_{[0,T]}) - \int_0^T D_s\left(-\sigma\frac{I_{(0)}}{I_{(1)}} + W_T\frac{I_{(0)}^2}{I_{(1)}^2} + \sigma\frac{I_{(0)}^2 I_{(2)}}{I_{(1)}^3}\right)ds \\ &= -\sigma W_T\frac{I_{(0)}}{I_{(1)}} + W_T^2\frac{I_{(0)}^2}{I_{(1)}^2} + \sigma W_T\frac{I_{(0)}^2 I_{(2)}}{I_{(1)}^3} + \int_0^T D_s\left(\sigma\frac{I_{(0)}}{I_{(1)}} - W_T\frac{I_{(0)}^2}{I_{(1)}^2} - \sigma\frac{I_{(0)}^2 I_{(2)}}{I_{(1)}^3}\right)ds.\end{aligned} \tag{6.4}$$

The Malliavin derivative of $W_T$ is $DW_T = \mathbb{1}_{[0,T]}$ $\mathbb{P}$-a.s. Thus it holds that

$$\int_0^T D_s W_T ds = T. \tag{6.5}$$

For the right-hand expression of the right-hand side of equation (6.4) we obtain

$$
\begin{aligned}
&\int_0^T D_s\left(\sigma\frac{I_{(0)}}{I_{(1)}} - W_T\frac{I_{(0)}^2}{I_{(1)}^2} - \sigma\frac{I_{(0)}^2 I_{(2)}}{I_{(1)}^3}\right)ds \\
&= \sigma\frac{I_{(1)}\sigma I_{(1)} - I_{(0)}\sigma I_{(2)}}{I_{(1)}^2} - \frac{I_{(1)}^2 T I_{(0)}^2 + 2I_{(1)}^2 W_T I_{(0)}\sigma I_{(1)} - 2W_T I_{(0)}^2 I_{(1)}\sigma I_{(2)}}{I_{(1)}^4} \\
&\quad - \sigma\frac{I_{(1)}^3 2I_{(0)}\sigma I_{(1)} I_{(2)} + I_{(1)}^3 I_{(0)}^2\sigma I_{(3)} - I_{(0)}^2 I_{(2)} 3I_{(1)}^2\sigma I_{(2)}}{I_{(1)}^6} \\
&= \sigma^2 - \sigma^2\frac{I_{(0)}I_{(2)}}{I_{(1)}^2} - T\frac{I_{(0)}^2}{I_{(1)}^2} - 2\sigma W_T\frac{I_{(0)}}{I_{(1)}} + 2\sigma W_T\frac{I_{(0)}^2 I_{(2)}}{I_{(1)}^3} - 2\sigma^2\frac{I_{(0)}I_{(2)}}{I_{(1)}^2} - \sigma^2\frac{I_{(0)}^2 I_{(3)}}{I_{(1)}^3} + 3\sigma^2\frac{I_{(0)}^2 I_{(2)}^2}{I_{(1)}^4} \\
&= \sigma^2 - 2\sigma W_T\frac{I_{(0)}}{I_{(1)}} - (TI_{(0)} + 3\sigma^2 I_{(2)})\frac{I_{(0)}}{I_{(1)}^2} + \sigma(2W_T I_{(2)} - \sigma I_{(3)})\frac{I_{(0)}^2}{I_{(1)}^3} + 3\sigma^2\frac{I_{(0)}^2 I_{(2)}^2}{I_{(1)}^4}.
\end{aligned}
\tag{6.6}
$$

From equation (6.6) it follows that the left-hand side of equation (6.4) can be written as

$$
\begin{aligned}
&-\sigma W_T\frac{I_{(0)}}{I_{(1)}} + W_T^2\frac{I_{(0)}^2}{I_{(1)}^2} + \sigma W_T\frac{I_{(0)}^2 I_{(2)}}{I_{(1)}^3} + \sigma^2 - 2\sigma W_T\frac{I_{(0)}}{I_{(1)}} - (TI_{(0)} + 3\sigma^2 I_{(2)})\frac{I_{(0)}}{I_{(1)}^2} \\
&+ \sigma(2W_T I_{(2)} - \sigma I_{(3)})\frac{I_{(0)}^2}{I_{(1)}^3} + 3\sigma^2\frac{I_{(0)}^2 I_{(2)}^2}{I_{(1)}^4} \\
&= \sigma^2 - 3\sigma W_T\frac{I_{(0)}}{I_{(1)}} + ((W_T^2 - T)I_{(0)} - 3\sigma^2 I_{(2)})\frac{I_{(0)}}{I_{(1)}^2} + \sigma(3W_T I_{(2)} - \sigma I_{(3)})\frac{I_{(0)}^2}{I_{(1)}^3} + 3\sigma^2\frac{I_{(0)}^2 I_{(2)}^2}{I_{(1)}^4}
\end{aligned}
\tag{6.7}
$$

As a consequence of equations (6.3), (6.4), and (6.7), we finally obtain

$$
\begin{aligned}
\pi_\Gamma &= \frac{1}{\sigma^2 x^2}\Bigg(-\sigma\pi_\Delta + \sigma^2 - 3\sigma W_T\frac{I_{(0)}}{I_{(1)}} + ((W_T^2 - T)I_{(0)} - 3\sigma^2 I_{(2)})\frac{I_{(0)}}{I_{(1)}^2} + \sigma(3W_T I_{(2)} - \sigma I_{(3)})\frac{I_{(0)}^2}{I_{(1)}^3} \\
&\quad + 3\sigma^2\frac{I_{(0)}^2 I_{(2)}^2}{I_{(1)}^4}\Bigg) \\
&= \frac{1}{\sigma^2 x^2}\Bigg(\sigma^2 - \sigma W_T\frac{I_{(0)}}{I_{(1)}} - \sigma^2\frac{I_{(0)}I_{(2)}}{I_{(1)}^2} + \sigma^2 - 3\sigma W_T\frac{I_{(0)}}{I_{(1)}} + ((W_T^2 - T)I_{(0)} - 3\sigma^2 I_{(2)})\frac{I_{(0)}}{I_{(1)}^2} \\
&\quad + \sigma(3W_T I_{(2)} - \sigma I_{(3)})\frac{I_{(0)}^2}{I_{(1)}^3} + 3\sigma^2\frac{I_{(0)}^2 I_{(2)}^2}{I_{(1)}^4}\Bigg) \\
&= \frac{1}{\sigma^2 x^2}\Bigg(2\sigma^2 - 4\sigma W_T\frac{I_{(0)}}{I_{(1)}} + ((W_T^2 - T)I_{(0)} - 4\sigma^2 I_{(2)})\frac{I_{(0)}}{I_{(1)}^2} + \sigma(3W_T I_{(2)} - \sigma I_{(3)})\frac{I_{(0)}^2}{I_{(1)}^3} \\
&\quad + 3\sigma^2\frac{I_{(0)}^2 I_{(2)}^2}{I_{(1)}^4}\Bigg).
\end{aligned}
$$

**Funding** Open Access funding enabled and organized by Projekt DEAL.

**Data Availability** All data used in our experiments have been produced with the R-packages Hudde (2021) and Gilbert and Varadhan (2015). The datasets generated and analysed during the current study are available from the corresponding author on reasonable request

## Declarations

**Conflict of Interests Statement** The authors Anselm Hudde and Ludger Rüschendorf certify that they have no affiliations with or involvement in any organization or entity with any financial interest (such as honoraria; educational grants; participation in speakers' bureaus; membership, employment, consultancies, stock ownership, or other equity interest; and expert testimony or patent-licensing arrangements), or non-financial interest (such as personal or professional relationships, affiliations, knowledge or beliefs) in the subject matter or materials discussed in this manuscript.

**Publisher's Note** Springer Nature remains neutral with regard to jurisdictional claims in published maps and institutional affiliations.